\documentclass[11pt,a4paper]{amsart}
\usepackage{array,amsmath,amssymb,graphicx}

\newtheorem{thm}{Theorem}
\newtheorem{lem}[thm]{Lemma}
\newtheorem{prop}[thm]{Proposition}

\newtheorem*{defi}{Definition}

\newcommand{\floor}[1]{\left \lfloor #1 \right \rfloor}

\begin{document}
\title{Crossing number of graphs and $\mathsf{\Delta Y}$-move}
\author{Youngsik Huh}
\address{Department of Mathematics,
College of Natural Sciences, Hanyang University, Seoul 04763,
Korea} \email{yshuh@hanyang.ac.kr}
\author{Ryo Nikkuni}
\address{Department of Mathematics, School of Arts and Sciences, Tokyo Woman’s Christian
University, 2-6-1 Zempukuji, Suginami-ku, Tokyo 167-8585, Japan}
\email{nick@lab.twcu.ac.jp}




\begin{abstract}
The crossing number of a graph is the minimum number of double points over all generic immersions of the graph into the plane. In this paper we investigate the behavior of crossing number under a graph transformation, called $\mathsf{\Delta Y}$-move, on the complete graph $K_n$. Concretely it is shown that for any $k\in \mathbb{N}$, there exist a natural number $n$ and a sequence of $\mathsf{\Delta Y}$-moves $K_n\rightarrow G^{(1)}\rightarrow \cdots \rightarrow G^{(k)}$ which is decreasing with respect to the crossing number. We also discuss the decrease of crossing number for relatively small $n$.
\end{abstract}

\maketitle



\section{Introduction}
When we consider graphs to be topological 1-complexes consisting of points (called {\em vertices}) and arcs connecting them (called {\em edges}) various topological notions of graphs are established. The crossing number and the genus are accepted as topological quantities measuring the nonplanarity of graphs~\cite{Sze}. The former comes from generic immersions of graphs into 2-dimensional spaces and the latter from embeddings into 2-dimensional spaces.

A {\em drawing} of a graph $G$ is a continuous map from $G$ to the Euclidean plane $\mathbb{R}^2$ (or the 2-sphere $\mathbb{S}^2$) such that its multiple points are only a finite number of transversal double points (called {\em crossings}) away from the image of vertices. Abusing the terminology, also the image of such a map will be called a {\em drawing}. Then the {\em crossing number} $\mathrm{cr}(G)$ of $G$ is defined to be the smallest number of crossings over all drawings of $G$. The determination of crossing number of graphs is hard in general~\cite{GJ, Hl, MC}.
The precise number is known for some specific families of graphs~\cite{EHK, Fi, BR, JS, Kle, Kle2}.
Even the famous two conjectures~\cite{Za,Guy} on the complete bipartite graphs $K_{p,q}$ and the complete graphs $K_n$ that
\begin{eqnarray*}
 & \mathrm{cr}(K_{p,q})=\floor{\frac{p}{2}}\floor{\frac{p-1}{2}}\floor{\frac{q}{2}}\floor{\frac{q-1}{2}}\;\mbox{and} \\
 & \mathrm{cr}(K_n)=\frac{1}{4} \floor{\frac{n}{2}} \floor{\frac{n-1}{2}} \floor{\frac{n-2}{2}} \floor{\frac{n-3}{2}}
\end{eqnarray*}
are proved only for $p\leq6$~\cite{Kl}, $(p,q)=(7,7), (7,8), (7,9)$~\cite{Woo} and $n\leq 12$~\cite{Guy2, PR}. Therefore it would be worthwhile to try to observe the behavior of crossing number under graph transformations.

As another sort of topological properties of graphs the intrinsic linkedness and the intrinsic knottedness can be taken. They are established from embeddings of graphs into the Euclidean 3-space $\mathbb{R}^3$~\cite{CG, Sa}.
A simple closed curve in $\mathbb{R}^3$ is called a {\em knot}. A knot is said to be {\em trivial} if it bounds a topological 2-dimensional disk. An $n$-component {\em link} is a disjoint union of $n$ simple closed curves which are simultaneously embedded in $\mathbb{R}^3$. A $2$-component link is said to be {\em splittable} if there exists a topological $2$-sphere in $\mathbb{R}^3$ which separates a component from the other. A graph is said to be {\em intrinsically linked} if every embedding of the graph into $\mathbb{R}^3$ contains a nonsplittable $2$-component link as its two disjoint cycles. An {\em intrinsically knotted} graph is a graph such that every its embedding into $\mathbb{R}^3$ contains a nontrivial knot as its cycle. It is known that $K_n$ is intrinsically linked and intrinsically knotted for $n\geq 6$ and $n\geq 7$, respectively~\cite{CG}. The intrinsically linked graphs were characterized in terms of graph minors~\cite{RST}.

\begin{figure}[t]
\includegraphics[scale=0.7]{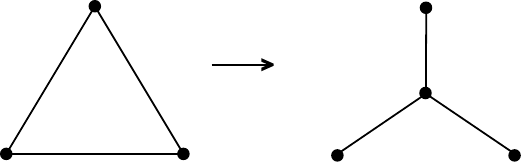}
\caption{$\mathsf{\Delta Y}$-move}
\label{fig1}
\end{figure}

A {\em $\mathsf{\Delta Y}$-move} is a transformation on graphs which replaces the three edges of a $3$-cycle by a $3$-star as depicted in Figure \ref{fig1}. Note that any $\mathsf{\Delta Y}$-move preserves the two properties: intrinsic linkedness and intrinsic knottedness.(Also the intrinsic linkedness is preserved by $\mathsf{Y\Delta}$-move.) On the other hand, a $\mathsf{\Delta Y}$-move increases the number of vertices and keeps the number of edges, hence it is expected that the crossing number decreases under the move. For example, if we perform $\mathsf{\Delta Y}$-moves on $K_7$ as many as possible then the resulting graph is the Heawood graph which is a cubic graph with 21 edges (See Figure \ref{fig3}). The crossing number of $K_7$ and the Heawood graph are $9$ and $3$\footnotemark, respectively.
In this paper, motivated by some observations like this, we investigate the behavior of crossing number under $\mathsf{\Delta Y}$-move.
\footnotetext{ A drawing of a graph can be considered to be a projected image onto $\mathbb{R}^2$ of an embedding of the graph into $\mathbb{R}^3$. It is known that every nontrivial knot produces at least three double points under the projection. Therefore, for every intrinsically knotted graph, $\mathrm{cr}\geq 3$.}
\begin{figure}[t]
\includegraphics[scale=0.6]{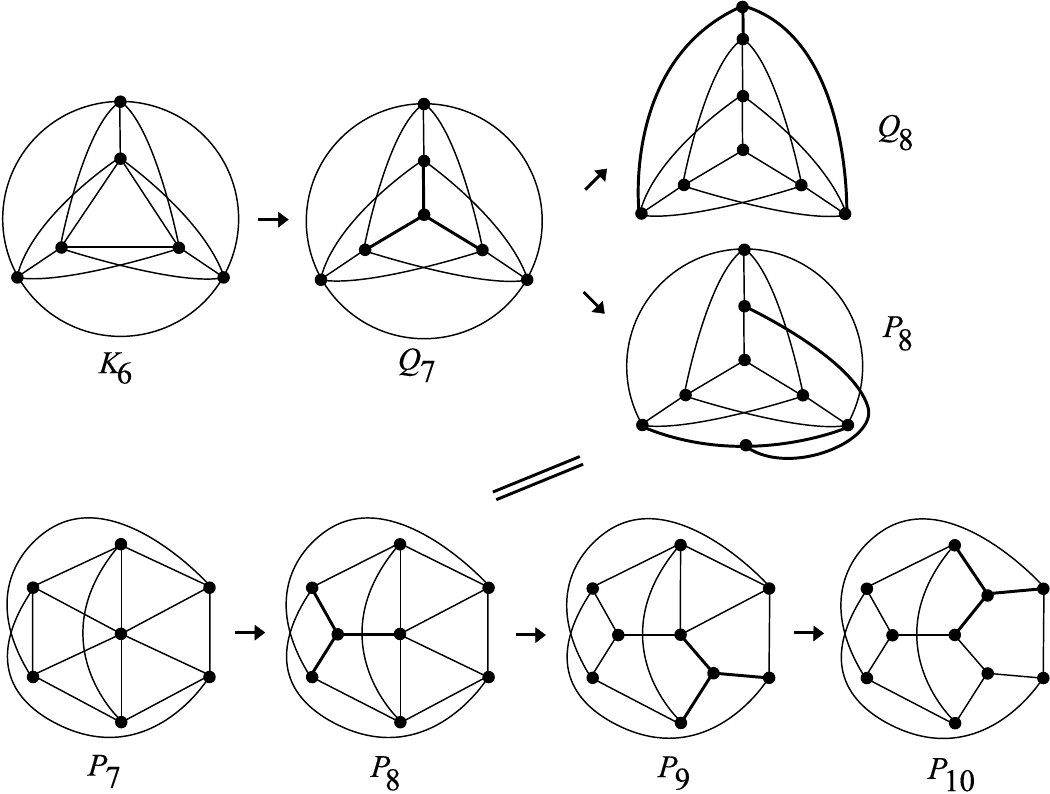}
\caption{Minimal-crossing drawings of  the Petersen Family: The family consists of  all graphs which are related to $K_6$ by $\mathsf{\Delta Y}$-moves.(The notations for the members follow \cite{HNTY}.)}
\label{fig2}
\end{figure}
\begin{figure}
\includegraphics[scale=0.6]{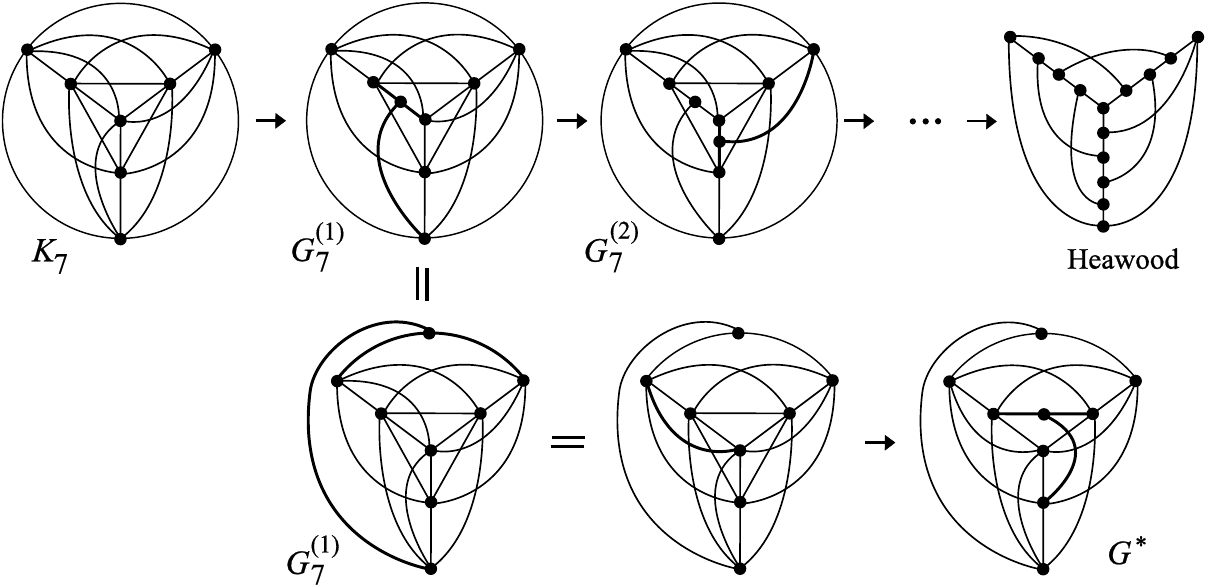}
\caption{Minimal-crossing drawings of $K_7$, $G_7^{(1)}$, $G_7^{(2)}$, $G^*$ and Heawood graph}
\label{fig3}
\end{figure}

\vspace{0.2cm}
Now the main results of this paper are described. Firstly we show
\begin{thm} \label{thm1}
For $n\geq 7$, a $\mathsf{\Delta Y}$-move on $K_n$ decreases the crossing number.
\end{thm}
\noindent Note that the theorem is not true for $n=6$. For the graph $Q_7$ in Figure \ref{fig2}, $\mathrm{cr}(Q_7)=3=\mathrm{cr}(K_6)$.

\vspace{0.2cm}
Let $G_n^{(1)}$ be the graph obtained from $K_n$ by a $\mathsf{\Delta Y}$-move on a $3$-cycle $\delta$ of $K_n$. Also let $\delta'$ be another $3$-cycle of $K_n$ which shares only one vertex with $\delta$, and $G_n^{(2)}$ be the graph obtained from $G_n^{(1)}$ by the $\mathsf{\Delta Y}$-move on $\delta'$. Then
\begin{thm} \label{thm2}
For $n\geq 7$, $\mathrm{cr}(G_n^{(2)}) < \mathrm{cr}(G_n^{(1)})$.
\end{thm}
\noindent Theorem \ref{thm2} may not be true when we select a $3$-cycle $\delta''$ (instead of $\delta'$) which is disjoint from $\delta$. In fact, for $n=7$,
$\mathrm{cr}(G_7^{(1)})=8=\mathrm{cr}(G^{*})$,
where $G^{*}$ is the graph obtained from $G_7^{(1)}$ by the $\mathsf{\Delta Y}$-move on $\delta''$ (See Figure \ref{fig3}).

\vspace{0.2cm}
Finally, as a generalization of Theorem \ref{thm1} and \ref{thm2}, we give
\begin{thm} \label{thm3}
For any $k\in \mathbb{N}$, there exist a natural number $n$ and a sequence of $\mathsf{\Delta Y}$-moves $K_n\rightarrow G^{(1)}\rightarrow \cdots \rightarrow G^{(k)}$ which is strictly decreasing with respect to the crossing number.
\end{thm}

\vspace{0.2cm}
The three theorems are proved in Section 3, 4 and 5, respectively. In Section 2 we introduce some necessary notions for our proofs.

The crossing number of the Petersen family given in Figure \ref{fig2}, although it is already known, can be shown by using the intrinsic linkedness. The proof is given in the final section for the readers' interest.
For the proof of $\mathrm{cr}(G_7^{(1)})=8=\mathrm{cr}(G^{*})$, the readers are referred to a computational verification system~\cite{CW}, because the authors' own proof is specific and not much simple.

\section{Crossing-reducible trigons} \label{sec2}
For a drawing of a graph, an {\em $n$-gon} will imply the image of an $n$-cycle of the graph on the drawing.
A drawing of a graph is said to be {\em good}, if it satisfies the three conditions in the below:
\begin{itemize}
\item[(G1)] No edge intersects itself, that is, there is no self-crossing.
\item[(G2)] There is no crossing between any two adjacent edges.
\item[(G3)] No two edges intersect each other more than once.
\end{itemize}
Note that if a drawing $D$ of a graph $G$ is a minimal-crossing drawing, that is, $\mathrm{cr}(D)=\mathrm{cr}(G)$, then it should be good. Also on a good drawing every trigon should be a simple closed curve.

\vspace{0.2cm}
Now we introduce a necessary notion for the proofs of the theorems.
Let $G$ be a graph, $\Delta_{v_1v_2v_3}$ be a 3-cycle of $G$, and $H$ be the graph obtained from $G$ by the $\mathsf{\Delta Y}$-move on $\Delta_{v_1v_2v_3}$. Here $v_1$, $v_2$ and $v_3$ are the vertices of the 3-cycle.

\begin{figure}
\includegraphics[scale=0.7]{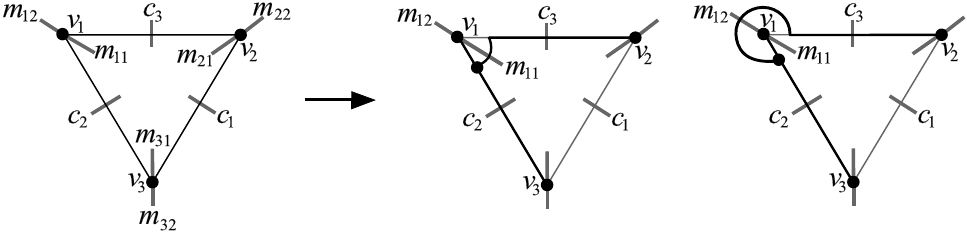}
\caption{ }
\label{fig4}
\end{figure}

Let $D$ be a good drawing of $G$.
We modify the drawing $D$ as depicted in Figure \ref{fig4}, and obtain drawings $D'$ and $D''$ of $H$.
The left-side of the figure illustrates a local picture of $D$ around the trigon $\Delta_{v_1v_2v_3}$, where $m_{11}$ and $m_{12}$ denote the numbers of  incident edges at the vertex $v_1$ which are locally going into the inside region and the outside region of the trigon, respectively.  For the edge $e(v_2,v_3)$ of $\Delta_{v_1v_2v_3}$ between $v_2$ and $v_3$, the number of crossings on $e(v_2,v_3)$ is denoted by $c_1$.

Now compare the crossing numbers of $D$, $D'$ and $D''$.
If we write
$$\mathrm{cr}(D)=c_1 + c_2 + c_3 + c^*\;,$$
where $c^*$ is the number of crossings which are not on any edge of $\Delta_{v_1v_2v_3}$, then
$$ \mathrm{cr}(D')=m_{11}+c_2+c_3+c^*\;\;\; \mbox{and} \;\;\;
\mathrm{cr}(D'')=m_{12}+c_2+c_3+c^*\;.$$
Therefore, if $c_1> \mbox{min}(m_{11}, m_{12})$, then
$$\mathrm{cr}(D)> \mbox{min}(\mathrm{cr}(D'),\mathrm{cr}(D'')) \geq \mathrm{cr}(H)\;.$$

Our observation can be summarized into the following lemma.
\begin{defi}
{\em A trigon $\Delta_{v_1v_2v_3}$ on a good drawing $D$ is a {\em cr-reducible trigon} of $D$, if
$c_i > \mathrm{min}(m_{i1}, m_{i2})$ for some $i$.}
\end{defi}
\begin{lem} \label{lem1}
Let $H$ be a graph obtained from a graph $G$ by the $\mathsf{\Delta Y}$-move on a 3-cycle $\Delta$. If $\Delta$ is a cr-reducible trigon of a good drawing $D$ of $G$,
then $\mathrm{cr}(D)>\mathrm{cr}(H)$.
\end{lem}

For our convenience we add a notation. For the trigon $\Delta_{v_1v_2v_3}$ in the above paragraph, let
$$\mathfrak{d}_{v_1v_2v_3}(v_i)=\mbox{min}(m_{i1}, m_{i2})\;.$$
Sometimes we denote the number by $\mathfrak{d}_{\Delta_{v_1v_2v_3}}(v_i)$.

\section{Proof of Theorem \ref{thm1}} \label{sec3}
Let $D$ be a minimal-crossing drawing of $K_n$.
We will prove Theorem \ref{thm1} by showing that $D$ has a cr-reducible trigon.

\vspace{0.2cm}
Firstly we consider the case $n=2k+1$ ($k\geq3$).
In this case, $\mathfrak{d}_{\Delta}(v) \leq k-1$ for every trigon $\Delta$ and every vertex $v$ of $\Delta$.

Suppose that $D$ has no cr-reducible trigon. Then, for every edge $e$,
$$\mathfrak{c}(e) \leq k-1\;, \;\;\;\;\cdots \cdots (*)$$
where $\mathfrak{c}(e)$ denotes the number of crossings of $D$ on the edge $e$.

\begin{figure}
\includegraphics[scale=0.8]{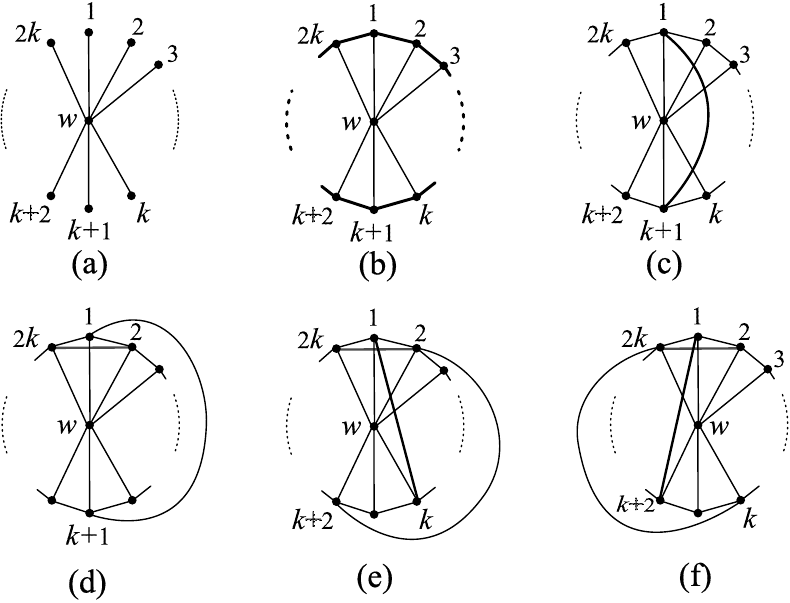}
\caption{$K_{2k+1}$ : Local pictures of $D$ around a vertex $w$. }
\label{fig5}
\end{figure}

Now select a vertex. By the two conditions (G1) and (G2) of good drawing,
the incident edges at the selected vertex constitute a spoke in $D$ as illustrated Figure \ref{fig5}-(a). Label the vertex by $w$, and the others by $1, \ldots ,2k$ (modulo $2k$) as in the figure. Then $\mathfrak{d}_{i,i+1,w}(w)=0$ in $D$, hence
$$\mathfrak{c}(e(i,i+1))=0\;\; \mbox{for every}\;\; i\;.$$
Therefore we see that any edge other than the edges drawn in \ref{fig5}-(b) should be  contained in the inside region or the outside of the $2k$-gon $P=P_{1,2,\ldots , 2k}$.

Suppose that $e(1,k+1)$ is contained in the inside of $P$. Then
$\mathfrak{c}(e(1,k+1))\geq k-1$, hence by ($*$),
$\mathfrak{c}(e(1,k+1))= k-1 $.
See Figure \ref{fig5}-(c). To avoid $e(1,k+1)$, the edges $e(2,k+1)$, $\ldots$, $e(2,2k)$ should be drawn in the outside of $P$.
This implies that the number of incident edges at {2} going into the inside of $\Delta_{2,3,w}$ is equal to $\mathfrak{d}_{2,3,w}(2)$. Note that all such edges should intersect $e(3,w)$.
To summarize,
$$\mathfrak{c}(e(3,w)) \geq \mathfrak{d}_{2,3,w}(2) +1\;,$$
that is, the trigon $\Delta_{2,3,w}$ is cr-reducible, which is a contradiction.

Now we can assume that every edge of the type $e(i,i+k)$ is contained in the outside of $P$. Then every two different $e(i,i+k)$ and $e(j,j+k)$ intersect each other, hence $\mathfrak{c}(e(i,i+k)) \geq k-1$. Again by ($*$),
$\mathfrak{c}(e(i,i+k)) = k-1$. Therefore each $e(i,i+k)$ has no more crossing with any edge of the other types. See Figure \ref{fig5}-(d), (e) and (f). If we draw $e(2,2k)$ and $e(1,k+1)$ together, then $e(2,2k)$ should be put into the inside of $P$. Also $e(1,k)$ and $e(1,k+2)$ should be there.
In conclusion $\mathfrak{c}(e(2,2k))\geq 3 > 1=\mathfrak{d}_{2,2k, w}(w)$, that is, $\Delta_{2,2k,w}$ is cr-reducible.

\begin{figure}
\includegraphics[scale=0.8]{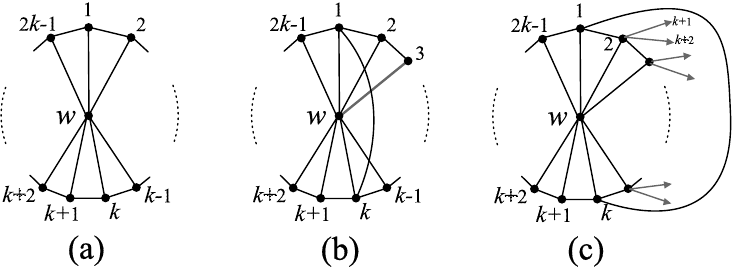}
\caption{$K_{2k}$ : Local pictures of $D$ around a vertex $w$. }
\label{fig6}
\end{figure}

\vspace{0.2cm}
Consider the case $n=2k$ ($k\geq 4$).
Suppose that $D$ has no cr-reducible trigon.
Then $\mathfrak{c}(e)\leq k-2$ for every edge $e$.
See Figure \ref{fig6}-(a).
Similarly with the previous case we see that
\begin{itemize}
\item[-] The incident edges at a vertex $w$ constitute a spoke in $D$.
\item[-] $\mathfrak{c}(e(i,i+1))=0$ for every $i$.
\item[-] Any other edge is contained in the inside region of the $(2k-1)$-gon $P$ or the outside region.
\end{itemize}

Suppose that $e(1,k)$ is contained in the inside region of $P$. See Figure \ref{fig6}-(b). Then $\mathfrak{c}(e(1,k)) \geq k-2$, hence
$\mathfrak{c}(e(1,k)) = k-2$, which implies that the edges $e(2,k+1)$, $\ldots$, $e(2, 2k-1)$ should be contained in the outside region of $P$.
Therefore $\mathfrak{c}(e(3,w)) \geq \mathfrak{d}_{2,3,w}(2) +1$, that is,
$\Delta_{2,3,w}$ is cr-reducible.

Now we can assume that every edge of the types $e(i,i+k-1)$ and $e(i,i+k)$ is contained in the outside of $P$. Then, as depicted in Figure \ref{fig6}-(c),
$$\mathfrak{c}(e(1,k)) \geq 2(k-2) > k-2 = \mathfrak{d}_{1,k,w}(w)\;,$$
which implies that $\Delta_{1,k,w}$ is cr-reducible.

\section{Proof of Theorem \ref{thm2}} \label{sec4}
\begin{figure}
\includegraphics[scale=0.8]{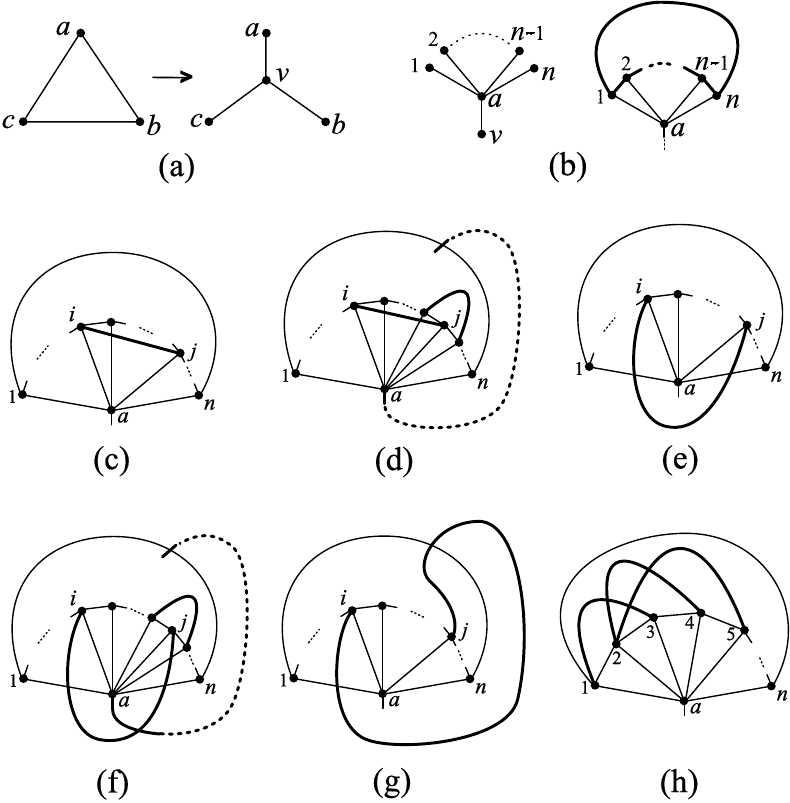}
\caption{ }
\label{fig7}
\end{figure}

For the complete graph $K_{n+3}$ with $n\geq 4$, let $\{a, b, c\}$ be the vertices of a 3-cycle of $K_{n+3}$ such that the $\mathsf{\Delta Y}$-move on the 3-cycle produces the graph $G_{n+3}^{(1)}$ as depicted in Figure \ref{fig7}-(a). The newly-born vertex of $G_{n+3}^{(1)}$ is denoted by $v$.
Let $D$ be a minimal-crossing drawing of $G_{n+3}^{(1)}$. We will prove Theorem \ref{thm2} by showing that $D$ has a cr-reducible trigon with one of $\{a, b, c\}$ as its vertex.

Suppose that $D$ has no cr-reducible trigon with one of $\{a, b, c\}$ as its vertex. We observe $D$ around the vertex $a$. Since $D$ is a good drawing, the incident edges at $a$ constitute a spoke in $D$. We label the vertices other than $\{v, a, b, c\}$ by $1$, $\ldots$, $n$ as illustrated in Figure \ref{fig7}-(b).
Then $\mathfrak{d}_{i,i+1,a}(a)=0$, hence

\vspace{0.2cm}
(\ref{sec4}-1) : $\mathfrak{c}(e(i,i+1))=0$ for every $1\leq i \leq n-1$.

\vspace{0.2cm}
\noindent Furthermore $\mathfrak{d}_{1,n,a}(a)=1$, hence

\vspace{0.2cm}
(\ref{sec4}-2) : $\mathfrak{c}(e(1,n)) \leq 1$.

\vspace{0.2cm}
From these two observations we can see that the $n$-gon $P_{1,2,\ldots , n}$ is a simple closed curve in $D$, and can assume that the edge $e(1,n)$ is drawn in $D$ as illustrated in Figure \ref{fig7}-(b).
Now we consider the position of an edge $e(i,j)$ in $D$ such that
$$j-i\geq 2 \;\; \mbox{and} \;\; \{i,j\} \neq \{1,n\}\;.$$

\vspace{0.2cm}
\noindent {\bf Case 1:} {\em $e(i,j)$ is contained in the inside region of the $(n+1)$-gon $P_{a,1,2, \ldots ,n}$. }

In this case, as depicted in Figure \ref{fig7}-(c)
$$\mathfrak{c}(e(i,j)) \geq j-i-1 \geq \mathrm{min}\{j-i-1, n-(j-i)\}=\mathfrak{d}_{i,j,a}(a)\;.$$
On the other hand, by the irreducibility of the trigon $\Delta_{ija}$,
$\mathfrak{c}(e(i,j))\leq \mathfrak{d}_{i,j,a}(a)$, hence it should be that
$$\mathfrak{c}(e(i,j)) = \mathfrak{d}_{i,j,a}(a) = j-i-1\;.$$
From (\ref{sec4}-1), (\ref{sec4}-2) and the equality in the above, we see that

\vspace{0.2cm}
{\em the vertex $v$ should be contained in the inside region of the $n$-gon $P_{1,2,  \cdots , n}$ and the edge $e(a,v)$ should intersect $e(1,n)$,}

\vspace{0.2cm}
\noindent because $v$ should be connected to the vertex $i+1$ along two edge-disjoint paths $v\rightarrow b \rightarrow i+1$ and  $v\rightarrow c \rightarrow i+1$.
Without loss of generality we assume $j\neq n$. Then the edge $e(j-1,j+1)$ should be contained in the inside region of $P_{123 \cdots n}$ as depicted in Figure \ref{fig7}-(d).

If the vertex $v$ is contained in the inside of $\Delta_{j-1,j,j+1}$, then the two paths $v\rightarrow b \rightarrow 1$ and  $v\rightarrow c \rightarrow 1$ intersect $e(j-1,j+1)$.
If $v$ is in the outside of $\Delta_{j-1,j,j+1}$, then $v\rightarrow b \rightarrow j$ and  $v\rightarrow c \rightarrow j$ intersect $e(j-1,j+1)$.
Therefore $\mathfrak{c}(e(j-1,j+1))\geq 2$, but $\mathfrak{d}_{j-1,j+1,a}(a)=1$. This contradicts the irreducibility of $\Delta_{j-1,j+1,a}$.

\vspace{0.2cm}
\noindent {\bf Case 2:} {\em $e(i,j)$ is contained in the outside region of the $n$-gon $P_{1,2, \cdots , n}$. }

We consider this case with Case 1 excluded. Then, by the goodness of $D$,  the edge $e(i,j)$ should intersect each of $e(1,a)$, $\ldots$, $e(i-1,a)$, $e(j+1,a)$, $\ldots$, $e(n,a)$ as depicted in Figure \ref{fig7}-(e).
To say again, $\mathfrak{c}(e(i,j)) \geq n-(j-i)-1$. Therefore, by the irreducibility of $\Delta_{i,j,a}$, it should be that

\vspace{0.2cm}
{\em $e(i,j)$ can have at most one more crossing other than these $n-(j-i)-1$ crossings.}

\vspace{0.2cm}
\noindent Now, without loss of generality, we assume $j\neq n$.

If $v$ is contained in the inside region of $P_{i,i+1, \ldots, j}$, then the two paths $v\rightarrow b \rightarrow n$ and  $v\rightarrow c \rightarrow n$ intersect $e(i,j)$, which contradicts our observation in the above.

If $v$ is contained in the outside region of $P_{i,i+1, \ldots, j}$ and the outside of $P_{1,2, \ldots, n}$, then the two paths $v\rightarrow b \rightarrow i+1$ and  $v\rightarrow c \rightarrow i+1$ intersect $e(1,n)$, which contradicts (\ref{sec4}-2) .

Lastly consider the case that $v$ is contained in the inside region of $P_{1,2, \ldots, n}$. Then $e(a,v)$ occupies the only possible additional crossing on $e(i,j)$, which implies that $e(j-1,j+1)$ should be contained in the inside region of $P_{1,2, \ldots, n}$ as depicted in Figure \ref{fig7}-(f). Repeating the same argument with Case 1, we can reach a contradiction.

\vspace{0.2cm}
We intend to conclude that every $e(i,j)$ is contained in the inside region of $P_{1,2,\ldots, n}$. So it needs to consider one more case.

\vspace{0.2cm}
\noindent {\bf Case 3:} {\em $e(i,j)$ intersects $e(1,n)$. }
(Note that this case happens only when $i\neq 1$ and $j \neq n$.)

See Figure \ref{fig7}-(g) for your understanding.
The only possible crossing on $e(1,n)$ is occupied by $e(i,j)$.
Therefore, excluding Case 1 and 2, all edges in the below should be contained in the inside region of $P_{1,2, \ldots, n}$.
$$\{e(1,k)\;|\;j+1 \leq k \leq n-1\} \;\; \mbox{and}  \;\;
\{e(l,k)\;|\; 2\leq l \leq j-1,\; j+1\leq k \leq n\}$$
Consequently,
$$\mathfrak{c}(e(i,j)) \geq 2+ (n-j-1) + (j-2)(n-j) > j-i-1 \geq \mathfrak{d}_{i,j,,a}(a)\;,$$
which contradicts the irreducibility of $\Delta_{ija}$.

\vspace{0.2cm}
Since all the three cases are excluded, we see that
{\em every $e(i,j)$ with $j-i\geq2$ and $\{i,j\}\neq\{1,n\}$ is contained in the inside region of $P_{1,2,\ldots ,n}$.}

\vspace{0.2cm}
For $n\geq5$, the above observation gives a local picture of $D$. As illustrated in Figure \ref{fig7}-(h)
the two edges $e(2,4)$ and $e(2,5)$ should intersect $e(1,3)$, hence
$$\mathfrak{c}(e(1,3))\geq 2 > 1= \mathfrak{d}_{13a}(a)\;,$$
which contradicts the irreducibility of $\Delta_{13a}$.

\vspace{0.4cm}
The remained case is that $n=4$. We will show that $\mathrm{cr}(D)\geq9$, which contradicts $\mathrm{cr}(G_7^{(1)})\leq8$ as seen in Figure \ref{fig3}.

\begin{figure}
\includegraphics[scale=1]{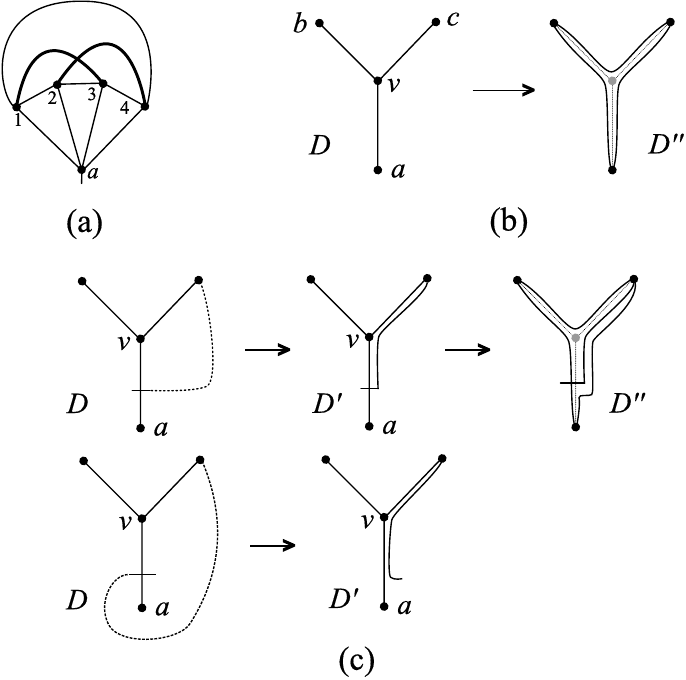}
\caption{ }
\label{fig8}
\end{figure}

Firstly we claim that $\mathfrak{c}(e(1,4))=0$.
By (\ref{sec4}-2), $\mathfrak{c}(e(1,4))\leq1$.
Since $e(1,3)$ and $e(2,4)$ should be contained in the inside region of $P_{1234}$, $e(1,3)$ and $e(2,4)$ intersect each other as illustrated in Figure \ref{fig8}-(a).
Then, by the irreducibility, $\mathfrak{c}(e(1,3))=1=\mathfrak{c}(e(2,4))$, that is, they have no other crossing.
But if $e(a,v)$ intersects $e(1,4)$, then the path $a \rightarrow v\rightarrow b\rightarrow 2$ intersects $e(1,3)$. If $e(b,v)$ (resp. $e(c,v)$) intersects $e(1,4)$, then
also $v\rightarrow b\rightarrow 2$ (resp. $v\rightarrow c\rightarrow 2$) intersects $e(1,3)$. If $e(b,2)$ intersects $e(1,4)$, then it intersects $e(1,3)$. It is same for $e(b,3)$, $e(c,2)$ and $e(c,3)$. Therefore there is no edge intersecting $e(1,4)$.

Define $a_{sp}$ to be
$$a_{sp} = \cup_{k=1}^4 e(a,k)\;,$$
and let $\mathfrak{c}(a_{sp}, b_{sp})$ denote the number of crossings of $D$
between $a_{sp}$ and $b_{sp}$.
Since  $\mathfrak{c}(e(1,3))=1=\mathfrak{c}(e(2,4))$ and
there is no crossing on $P_{1234}$, we can see
$$\mathfrak{c}(a_{sp},b_{sp}) + \mathfrak{c}(b_{sp},c_{sp}) + \mathfrak{c}(c_{sp},a_{sp}) \geq 2 + 2 + 2 \geq 6\;.$$

Now suppose that $\mathrm{cr}(D)\leq 8$. Then
$$\mathfrak{c}(a_{sp},b_{sp}) + \mathfrak{c}(b_{sp},c_{sp}) + \mathfrak{c}(c_{sp},a_{sp}) \;=\;7 \;\; \mbox{or}\;\; 6\;.$$
In the former case, it should be that $\mathrm{cr}(D)=8$. Hence $\mathfrak{c}(e(a,v))=0=\mathfrak{c}(e(b,v))=\mathfrak{c}(e(c,v)),$
which enables us to get a drawing $D''$ of $K_7$ with $\mathrm{cr}(D'')=8$ as illustrated in Figure \ref{fig8}-(b).

In the latter case, if $\mathrm{cr}(D)=7$, then we also have
$\mathfrak{c}(e(a,v))=0=\mathfrak{c}(e(b,v))=\mathfrak{c}(e(c,v))$,
which gives the same contradiction with the former case.
When $\mathrm{cr}(D)=8$, without loss of generality, we can assume that
$\mathfrak{c}(e(a,v))=1$ and $\mathfrak{c}(e(b,v))=\mathfrak{c}(e(c,v))=0$.
The edge intersecting $e(a,v)$ at the crossing should be incident to $b$ or $c$, because there is no crossing on $P_{1234}$. Then, by modifying $D$ as illustrated in Figure \ref{fig8}-(b) and (c), we can obtain a drawing of $D'$ of $G_7^{(1)}$ or $D''$ of $K_7$ such that $\mathrm{cr}(D')<\mathrm{cr}(D)$ or $\mathrm{cr}(D'')\leq 8$.

\section{Proof of Theorem \ref{thm3}} \label{sec5}
\begin{figure}
\includegraphics[scale=0.9]{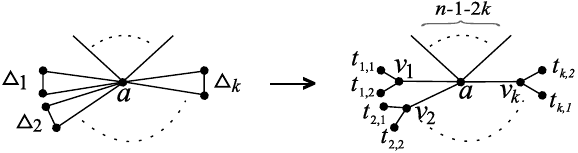}
\caption{$G_n^{(k)} $is obtained from $K_n$ by $\mathsf{\Delta Y}$-moves on $\Delta_1$, $\ldots$, $\Delta_k$.  }
\label{fig9}
\end{figure}

Select a vertex $a$ of $K_n$.
Let $\Delta_1$, $\ldots$, $\Delta_k$ $(k\geq 2)$ be mutually edge-disjoint 3-cycles of $K_n$ sharing the vertex $a$. Then $G_n^{(k)}$ will denote the graph obtained from $K_n$ by $\mathsf{\Delta Y}$-moves on all of the 3-cycles $\Delta_1$, $\ldots$, $\Delta_k$. As depicted in Figure \ref{fig9}, each $v_i$ denotes the vertex of $G_n^{(k)}$
which is born from one of the $\mathsf{\Delta Y}$-moves.
The other vertices adjacent to $a$ will be  labelled simply by the integers $\{1,\ldots , n-1-2k\}$.
The two vertices other than $a$ which are adjacent to $v_i$ are denoted by $t_{i,1}$ and $t_{i,2}$.

Note that Theorem \ref{thm2} is identical with Theorem \ref{thm3} for $k=2$.
We will prove Theorem \ref{thm3} by showing the following statement:
\begin{center}
{\em For $k \geq 2$ and $n \geq 10k+1$, $\mathrm{cr}(G_n^{(k)}) > \mathrm{cr}(G_n^{(k+1)})$.}
\end{center}

\begin{figure}
\includegraphics[scale=0.9]{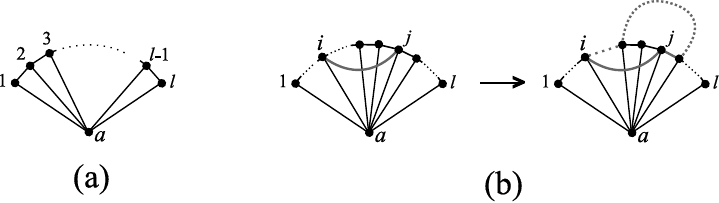}
\caption{ }
\label{fig10}
\end{figure}

\vspace{0.2cm}
Let  $D$ be a minimal-crossing drawing of $G_n^{(k)}$.
Supposing that $D$ has no cr-reducible trigon with $a$ as its vertex, 
we pursue a contradiction.
The edges $e(a,1)$, $\ldots$, $e(a,n-1-2k)$ will be called {\em raw edges} at $a$. Then, in the drawing $D$, for some $l\geq \frac{n-1-2k}{k}$, we find $l$ raw edges which are consecutive around $a$. The vertices of the raw edges other than $a$ are labelled by $1$, $\ldots$, $l$ as illustrated in Figure \ref{fig10}-(a). Since $\mathfrak{d}_{i,i+1,a}(a)=0$ in $D$ for each $1 \leq i \leq l-1$, the cr-irreducibiliy implies
$$\mathfrak{c}(e(i,i+1))=0\;.$$

Now we consider the position of an edge $e(i,j)$ such that $1<i<j<l$ and $j-i\geq3$. For two edge-disjoint subgraphs $H$and $H'$ of a graph $G$, let $\mathfrak{c}(H,H')$ denote the number of crossings of $D$ between $H$ and $H'$.

\vspace{0.2cm}
\noindent {\bf Case 1:} {\em $e(i,j)$ is contained in the inside region of the polygon $P_{a, 1, 2, \ldots , l}$. }

In this case, $\mathfrak{c}(e(i,j))\geq j-i-1$. Hence, by the irreducibility and the goodness of $D$,
$$\mathfrak{c}(e(i,j))= j-i-1=\mathfrak{d}_{ija}(a)\;,$$
as illustrated in Figure \ref{fig10}-(b). Also, by the irreducibility,
$$\mathfrak{c}(e(j-2, j+1)) \leq 2\;.$$

Now we observe the crossings on the polygon 
$$P_*=i \rightarrow i+1 \rightarrow \cdots \rightarrow j-2 \rightarrow j+1 \rightarrow j \rightarrow i \;.$$
The path $i\rightarrow i+1 \rightarrow \cdots \rightarrow j-2$ has no crossing on itself, and $\mathfrak{c}(e(j+1,j))=0$. The edge $e(j,i)$ is intersected only by $e(a,i+1)$, $e(a,i+2)$, $\ldots$, $e(a,j-1)$. Therefore any more crossings on $P_*$ should be positioned only on $e(j-2,j+1)$.

From the observation in the above we claim that $P_*$ separates the vertex $j-1$ from $a$. Otherwise $\mathfrak{c}(P_*, e(a,j-1))$ should be an even number. Then, since
$$\mathfrak{c}(P_*, e(a,j-1)) \geq \mathfrak{c}(e(j,i),e(a,j-1))=1\;,$$
$e(j-2,j+1)$ should intersect $e(a,j-1)$, hence also intersect  $e(a,j-2)$ and $e(a,j)$. Consequently $\mathfrak{c}(e(j-2,j+1))\geq 3$, which is a contradiction. (In fact it also contradicts the goodness of $D$.)

Since $P_*$ separates the vertex $j-1$ from $a$,
each path $a\rightarrow v_{m_1} \rightarrow t_{m_1,m_2} \rightarrow j-1$ ($1\leq m_1 \leq k$, $m_2=1, 2$) should intersect $e(j-2, j+1)$. Therefore, for $k\geq 3$, we have the contradiction $\mathfrak{c}(e(j-2,j+1))\geq 3$.
For $k=2$, the paths $a\rightarrow v_{m_1} \rightarrow t_{m_1,m_2} \rightarrow j-1$ contribute at least 2 to $\mathfrak{c}(e(j-2,j+1)$. Since the even number $\mathfrak{c}(P_*, \Delta_{1, j-1,a})$ is not zero, we have $\mathfrak{c}(e(j-2,j+1))\geq 3$ again.

\vspace{0.2cm}
\noindent {\bf Case 2:} {\em $e(i,j)$ intersects the inside region of the polygon $P_{a, 1, 2, \cdots , l}$. }

\begin{figure}
\includegraphics[scale=0.7]{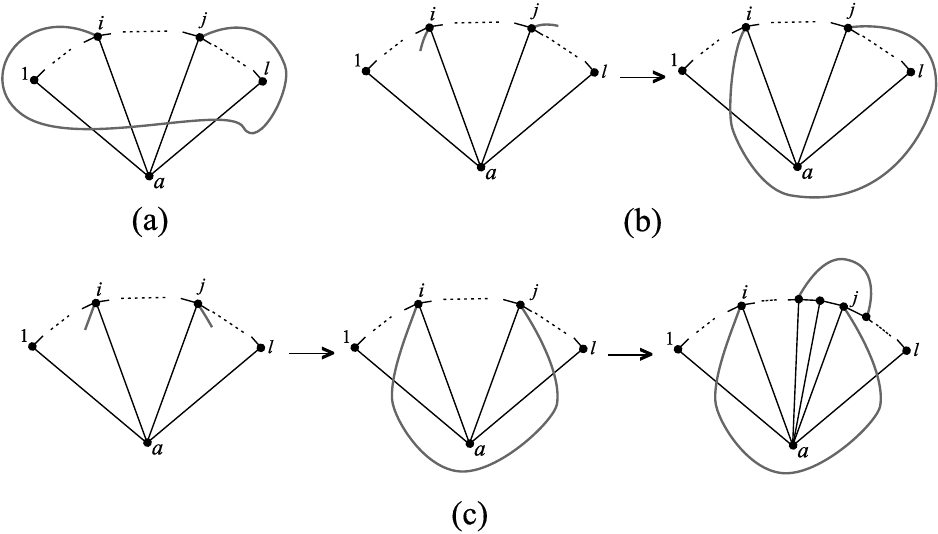}
\caption{ }
\label{fig11}
\end{figure}

\vspace{0.2cm}
Note that we consider Case 2 to reach the conclusion:  {\em every $e(i,j)$ is contained in the outside region of $P_{a, 1, 2, \ldots ,l}$.} For the purpose, under the exclusion of Case 1, we divide Case 2 into the three subcases which are depicted in Figure \ref{fig11}-(a):
\begin{itemize}
\item[(S1)]{\em $e(i,j)$ is locally contained in the outside region near $i$ and $j$.}
\item[(S2)]{\em $e(i,j)$ is locally contained in the inside region near $i$ and  in the outside region near $j$.}
\item[(S3)]{\em $e(i,j)$ is locally contained in the inside region near $i$ and $j$.}
\end{itemize}
For our discussion recall the inequality from the cr-irreducbility
$$\mathfrak{c}(e(i,j))\leq \mathfrak{d}_{ija}(a)= \mbox{min}\{j-i-1, n-(j-i)-k-2\} < l \;\;.$$

\vspace{0.2cm}
The subcase (S1) cannot happen, because $\mathfrak{c}(e(i,j))\geq l$ by the goodness (G3) when $e(i,j)$ intersects the inside region.

Consider the subcase (S2). By the inequality in the above and (G3),
the edge $e(i,j)$ intersects only the edges
$e(a,1)$, $e(a,2)$, $\ldots$, $e(a,i-1)$ among the edges
$e(a,1)$, $e(a,2)$, $\ldots$, $e(a,l)$
as depicted in Figure \ref{fig11}-(b). Therefore the polygon $P_{i,i+1, \ldots, j}$ separates $\{1, \ldots, i-1\}$ from $\{j+1, \ldots, l\}$ and $\{a\}$.
Furthermore, for each $l+1 \leq m_1 \leq n-1-2k$, the set of paths
$$\{a\rightarrow m_1 \rightarrow m_2 \;|\; 1\leq m_2 \leq  i-1\}$$
contributes at least one crossing on $e(i,j)$. Therefore we have
$$\mathfrak{c}(e(i,j)) \geq  (i-1)(l-j)+ (i-1)+ (n-1-2k-l)\;.$$
There are more crossings on $e(i,j)$. If some $v_{m_3}$ is in the inside region of $P_{i,i+1, \ldots, j}$, then the two paths $a\rightarrow v_{m_3} \rightarrow t_{m_3,1} \rightarrow 1$ and  $a\rightarrow v_{m_3} \rightarrow t_{m_3,2} \rightarrow 1$ intersect $e(i,j)$. If $v_{m_3}$ is in the outside region, then  $a\rightarrow v_{m_3} \rightarrow t_{m_3,1} \rightarrow l$ and  $a\rightarrow v_{m_3} \rightarrow t_{m_3,2} \rightarrow l$ intersect $e(i,j)$. In consequence
$$\mathfrak{c}(e(i,j)) \geq (i-1)(l-j+1)+ (n-1-2k-l) + 2k \;.$$
Since the rightside in the inequality is greater than $n-(j-i)-k-2$, we have a contradiction
$$\mathfrak{c}(e(i,j)) >  n-(j-i)-k-2 \geq \mathfrak{d}_{ija}(a)\;.$$

Figure \ref{fig11}-(c) depicts the situation of subcase (S3).
Since $e(i,j)$ can not intersect $e(a,i)$ and $e(a,j)$,
it should intersect $e(a,1)$, $\ldots$, $e(a,i-1)$, $e(a,j+1)$. $\ldots$, $e(a,l)$.
Therefore the polygon $P_{i,i+1, \dots , j}$ separates the vertex $a$ from
$\{1, \ldots , i-1\}$ and $\{j-1, \ldots , l \}$, hence the paths $a\rightarrow l+1 \rightarrow 1$, \ldots $a\rightarrow n-1-2k \rightarrow 1$ intersect $e(i,j)$. In consequence
$$\mathfrak{c}(e(i,j)) \geq (i-1)+(l-j)+(n-1-2k-l) = n-(j-i)-2k-2 \;.$$
Again, by the irreduciblity of $\Delta_{a,i,j}$,
$$\mathfrak{c}(e(i,j))  \leq \mathfrak{d}(a)_{a,i,j} \leq n-(j-i)-k-2 \;.$$
Therefore $e(i,j)$ can have at most $k$ crossings more on itself.
In fact $\mathfrak{c}(e(i,j))  =n-(j-i)-k-2 $, because
the additional $k$ crossings come from the intersection with
$$\{a\rightarrow v_m \rightarrow t_{m,1}\rightarrow 1, \;a\rightarrow v_m \rightarrow t_{m,2}\rightarrow 1 \},\;\;\;
1\leq m \leq k \;.$$
If some $v_m$ is contained in the inside region of $P_{i,i+1, \ldots , j}$, then
both the two edge-disjoint paths $v_m \rightarrow t_{m,1} \rightarrow 1$ and
$v_m \rightarrow t_{m,2} \rightarrow 1$ intersect $e(i,j)$ which is contradictory to
$\mathfrak{c}(e(i,j))  =n-(j-i)-k-2 $.
In conclusion $\{v_m, t_{m,1}, t_{m,2} \}$ should be contained in the outside region of
$P_{i,i+1, \ldots , j }$ for every $m$.

From the counting of crossings on $e(i,j)$, we know that the edge $e(j-2, j+1)$ can not intersect $e(i,j)$, hence $e(j-2, j+1)$ is located in the outside regions of $P_{a,1,2, \ldots , l}$ and $P_{i,i+1, \ldots , j}$ as illustrated in the last figure of Figure \ref{fig11}-(c). From the intersections with $e(1,j-1)$, $e(l,j-1)$ and 
$$\{a \rightarrow v_m \rightarrow t_{m,m'} \rightarrow j-1 \; |  1\leq m \leq k\} \;,$$ 
it should be that $\mathfrak{c}(j-2,j+1) > 2$. On the contrary, by the irreducibility, 
$$\mathfrak{c}(j-2,j+1) \leq \mathfrak{d}(a)_{a, j-2, j+1} \leq 2\;.$$

\begin{figure}
\includegraphics[scale=0.9]{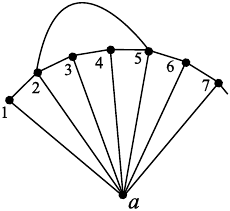}
\caption{ }
\label{fig12}
\end{figure}

\vspace{0.2cm}
Since Case 1 and 2 are excluded, it can be assumed that every $e(i,j)$ with $1<i<j<l$ and $j-i \geq 3$ is located in the outside region of $P_{a,1,2, \ldots , l}$. Now we set $l\geq 8$.
Then $\mathfrak{c}(e(2,5))$ intersects $e(3,6)$, $e(3,7)$ and $e(4,7)$ as illustrated in Figure \ref{fig12}, which contradicts the irreduciblity $\mathfrak{c}(e(2,5)) \leq \mathfrak{d}_{a,2,5}(a) \leq 2$.

\vspace{0.2cm}
\noindent \textbf{Remark.} In the proof we set $l\geq \frac{n-1-2k}{k}$ and $l\geq 8$. Therefore if $n\geq 10k+1$, then the existence of crossing-deceasing $\mathsf{\Delta Y}$-moves with length $k$ is guaranteed. But, as you can see from Theroem \ref{thm1} and \ref{thm2}, $n=7$ is enough for $k=2$.
As a future study we may try to find more optimal $n$.

\section{The intrinsic linkedness and crossing number of the Petersen family}
The Petersen family consists of all graphs which are related to $K_6$ by $\mathsf{\Delta Y}$ and $\mathsf{Y \Delta}$-moves. Figure \ref{fig2} shows their drawings.
Since the graphs are relatively small there would be several ways to determine their crossing numbers.
In this section, for the readers' interest, we determine the crossing number of the Petersen family graphs by using the intrinsic linkedness.
\begin{prop}\label{prop1} 
$\mbox{ }$ \\
(1) For $P_8$, $P_9$ and $P_{10}$, $\mathrm{cr}=2$. \\
(2) For $K_6$, $Q_7$, $P_7$ and $Q_8$, $\mathrm{cr}=3$.
\end{prop}

\begin{figure}
\includegraphics[scale=0.6]{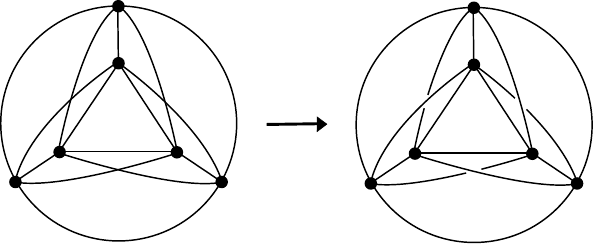}
\caption{A drawing $f(G)$ and a lift $\hat{f}(G)$}
\label{fig13}
\end{figure}

For a drawing $f:G \rightarrow \mathbb{R}^2$ of a graph $G$, by adding information on which strand passes over/under at each crossing, we obtain a diagram which represents an embedding $\hat{f}: G \rightarrow \mathbb{R}^3$, called a {\em lift} of $f$ (See Figure \ref{fig13}). Note that $\pi \circ \hat{f}=f$ for the natural projection $\pi:\mathbb{R}^3 \rightarrow \mathbb{R}^2$.

\vspace{0.2cm}
Now we prove Proposition \ref{prop1}. Let $G$ be a graph in the Petersen family and $f$ be a drawing of $G$. Note that $K_6$ is intrinsically linked and the intrinsic linkedness is preserved by  $\mathsf{\Delta Y}$ and $\mathsf{Y\Delta}$-moves \cite{RST, Sa}. Hence every lift of $f$ contains a nonsplittable $2$-component link. It is known that the projected image of any nonsplittable link under $\pi$ contains at least two crossings, which implies $\mathrm{cr}(f) \geq 2$. Therefore (1) of Proposition \ref{prop1} comes from the drawings in Figure \ref{fig2}.

\begin{figure}
\includegraphics[scale=0.9]{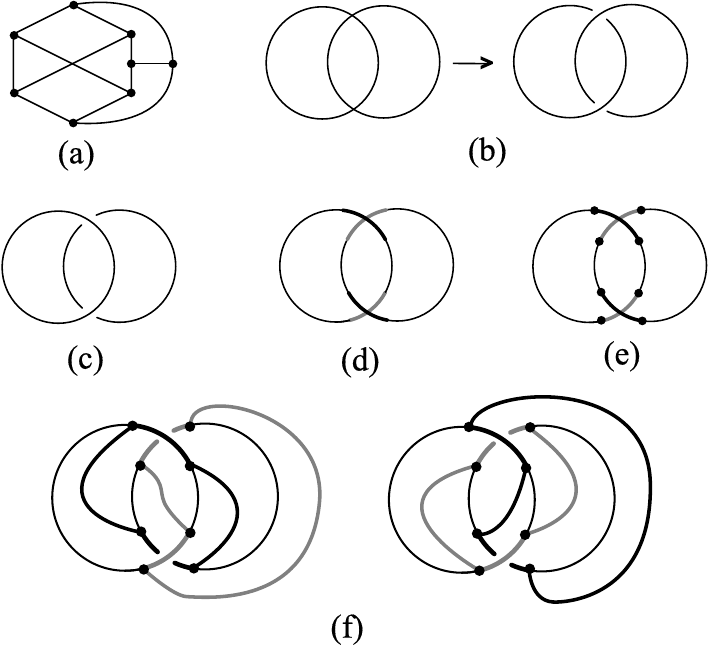}
\caption{(a): The graph $H$. (b): $f(C_1 \cup C_2)$ and $\hat{f}(C_1 \cup C_2)$. (c): $\tilde{f}(C_1 \cup C_2)$. (d): $D_1 \cup D_2$ share the four subarcs with $C_1 \cup C_2$ near the crossing points. (e): Each $C_i$ has four vertices near the crossing points. (f): Two possible shapes of $\tilde{f}(C_1 \cup C_2 \cup D_1 \cup D_2)$.}
\label{fig14}
\end{figure}

For the proof of (2) of Proposition \ref{prop1} we introduce a lemma.
\begin{lem} \label{lem2}
For an intrinsically linked graph $G$ with $\mathrm{cr}(G)=2$, \\
(1) $|V(G)|\geq 8$. \\
(2) Furthermore, if $|V(G)|= 8$, then $G$ contains a subgraph $H$ which is illustrated in Figure \ref{fig14}-(a).
\end{lem}
From Figure \ref{fig2} we see that $\mathrm{cr} \leq 3$ for $K_6$, $Q_7$, $Q_8$ and $Q_8$. On the other hand the first statement of Lemma \ref{lem2} implies that $\mathrm{cr} \geq 3$ for $K_6$, $Q_7$ and $P_7$.
Lastly, the graph $H$ is not a subgraph of $Q_8$, because $H$ contains a cycle of length $5$ but $Q_8$ does not.
Therefore, by (2) of Lemma \ref{lem2}, $\mathrm{cr}(Q_8) \geq 3$.

\vspace{0.2cm}
\noindent {\em Proof of Lemma \ref{lem2}.}
Let $f$ be a drawing of $G$ with $\mathrm{cr}(f)=2$, and $\hat{f}$ be a lift of $f$. Since $G$ is intrinsically linked and $\mathrm{cr}(f)=2$, the graph has two disjoint cycles $C_1$ and $C_2$ such that $f(C_1 \cup C_2)$ and $\hat{f}(C_1 \cup C_2)$ are the shapes illustrated in Figure \ref{fig14}-(b).
By interchanging the over and under strands at a crossing of $\hat{f}(C_1 \cup C_2)$ we obtain another lift $\tilde{f}$ such that $\tilde{f}(C_1 \cup C_2)$ is splittable as illustrated in Figure \ref{fig14}-(c).
By the intrinsic linkedness of $G$ we see that $G$ should have another pair of disjoint cycles $D_1$ and $D_2$ such that $\tilde{f}(D_1 \cup D_2)$ is nonsplittable again. The drawing $f$ has only two crossings, hence $f(D_1 \cup D_2)$ share the four subarcs of $f(C_1 \cup C_2)$ which are illustrated as thick curves in Figure \ref{fig14}-(d). Furthermore, for $D_1$ and $D_2$ to be disjoint, each of $C_1$ and $C_2$ contains at least four vertices as illustrated in Figure \ref{fig14}-(e). Therefore $|V(G)|\geq 8$.

Finally assume that $|V(G)|= 8$, and try to add $f(D_1 \cup D_2)$ onto $f(C_1 \cup C_2)$ without producing any more crossing and vertex. Then $\tilde{f}(C_1\cup C_2 \cup D_1 \cup D_2)$ should be one of the shapes in Figure \ref{fig14}-(f), from which we see that $C_1\cup C_2 \cup D_1 \cup D_2$ is isomorphic to $H$.

\section*{Acknowledgements}
The first author was supported by the National Research Foundation of Korea (NRF) grant funded by the Korea government (MSIP) (NRF - 2016R1D1A1B01008044).

The second author was supported by JSPS KAKENHI Grant Number JP22K03297.




\begin{thebibliography}{XX}
\bibitem{BR} L. W. Beineke and R. D. Ringeisen, {\em On the crossing numbers of products of cycles and graphs of order four}, J. Graph Theory \textbf{4} (1980), 145–-155.
\bibitem{CG} J. H. Conway and C. McA. Gordon, {\em Knots and links in spatial graphs}, J. Graph Theory  \textbf{7} (1983), 445–-453.
\bibitem{CW} M. Chimani and T. Wiedera, {\em An ILP-based proof system for the crossing number problem}, In Proc. ESA 2016, volume 57 of LIPIcs, pages 29:1--29:13, 2016.
\bibitem{EHK} G. Exoo, F. Haray and J. Kabell, {\em The crossing numbers of some generalized Petersen graphs} Math. Scand. \textbf{48} (1981), 184–-188.
\bibitem{Fi} S. Fioroni, {\em On the crossing number of generalized Petersen graphs}, Combinatorics '84 (Bari, 1984), 225–-241. North-Holland Math. Stud., 123 Ann. Discrete Math., 30 North-Holland Publishing Co., Amsterdam, 1986.
\bibitem{GJ} M.R. Garey and D. S. Johnson, {\em Crossing number is NP-complete},
SIAM J. Algebraic Discrete Methods \textbf{4} (1983), 312--316.
\bibitem{Guy} R. K. Guy, {\em A combinatorial problem}, Nabla (Bulletin of the Malayan Mathematical Society) \textbf{7} (1960), 68–-72.
\bibitem{Guy2} R. K. Guy, {\em Crossing numbers of graphs}, In {\em Graph theory and applications (Proc. Conf., Western Michigan Univ., Kalamazoo, Mich., 1972; dedicated to the memory of J. W. T. Youngs)}, 111–-124. Lecture Notes in Math., Vol. 303 Springer-Verlag, Berlin-New York, 1972.
\bibitem{Hl} P. Hlin\v{e}n\'{y}, {\em Crossing number is hard for cubic graphs},
J. Combin. Theory Ser. B \textbf{96} (2006), 455–-471.
\bibitem{HNTY} R. Hanaki, R. Nikkuni, K. Taniyama and A. Yamazaki,
{\em On intrinsically knotted or completely 3-linked graphs},
Pacific J. Math. 252 (2011), 407--425.
\bibitem{JS} S. Jendrol' and M. Ščerbová, {\em On the crossing numbers of $S_m \times P_n$ and $S_m \times C_n$}, Časopis Pěst. Mat. \textbf{107} (1982), 225–-230.
\bibitem{Kl} D. J. Kleitman, {\em The crossing number of $K_{5,n}$},
J. Combinatorial Theory \textbf{9} (1970), 315–-323.
\bibitem{Kle} M. Klešč, {\em The crossing numbers of products of paths and stars with 4-vertex graphs}, J. Graph Theory \textbf{18} (1994), 605–-614.
\bibitem{Kle2} M. Klešč, {\em The crossing numbers of Cartesian products of paths with 5-vertex graphs}, Graph theory (Prague, 1998) Discrete Math. \textbf{233} (2001), 353-–359.
\bibitem{MC} S. Cabello and B. Mohar, {\em Adding one edge to planar graphs makes crossing number and 1-planarity hard}, SIAM J. Comput. \textbf{42} (2013), 1803–-1829.
\bibitem{PR} S. Pan and R. B. Ritchter, {\em The crossing number of $K_{11}$ is $100$},
J. Graph Theory \textbf{56} (2007), 128–-134.
\bibitem{RST} N. Robertson, P. Seymour and R. Thomas, {\em Sachs’ linkless embedding conjecture}, J.
Combin. Theory Ser. B \textbf{64} (1995), 185–-227.
\bibitem{Sa} H. Sachs,
{\em On spatial representations of finite graphs}, 
Finite and infinite sets, Vol. I, II (Eger, 1981), 649–-662.
Colloq. Math. Soc. János Bolyai, 37
North-Holland Publishing Co., Amsterdam, 1984
\bibitem{Sze} L. A. Sz\'ekely, {\em A successful concept for measuring non-planarity of graphs: the crossing number}, 6th International Conference on Graph Theory
Discrete Math. \textbf{276} (2004), no.1--3, 331–-352.
\bibitem{Woo} D. R. Woodall, 
{\em Cyclic-order graphs and Zarankiewicz's crossing-number conjecture},
J. Graph Theory \textbf{17} (1993) 657–-671.
\bibitem{Za} K. Zarankiewicz, {\em On a problem of P. Turan concerning graphs}, Fund. Math. \textbf{41} (1954), 137–-145.


\end{thebibliography}
\end{document}